\documentclass{amsart}
\usepackage{amssymb,latexsym,amsthm,amsmath}
\usepackage{commath}

\begin{document}

\newtheorem{theorem}{Theorem}[section]
\newtheorem{lemma}[theorem]{Lemma}
\newtheorem{proposition}[theorem]{Proposition}

\theoremstyle{definition}
\newtheorem{definition}[theorem]{Definition}

\theoremstyle{remark}
\newtheorem{remark}[theorem]{Remark}

\numberwithin{equation}{section}

\newcommand{\C}{\mathbb{C}}
\newcommand{\T}{\mathbb{T}}
\newcommand{\PnC}{\mathbb{P}^n(\C)}
\newcommand{\N}{\mathbb{N}}
\newcommand{\Z}{\mathbb{Z}}
\newcommand{\mA}{\mathcal{A}}
\newcommand{\mP}{\mathcal{P}}
\newcommand{\mT}{\mathcal{T}}
\newcommand{\U}{\mathrm{U}}
\newcommand{\PU}{\mathrm{PU}}
\newcommand{\SU}{\mathrm{SU}}
\newcommand{\GL}{\mathrm{GL}}
\newcommand{\SL}{\mathrm{SL}}
\newcommand{\Spe}{\mathrm{S}}
\newcommand{\End}{\mathrm{End}}

\title[Toeplitz Operators and Representation Theory]{Separately Radial and Radial Toeplitz Operators on the Projective Space and Representation Theory}
\author{R.~Quiroga-Barranco}
\address{Centro de Investigaci\'on en Matem\'aticas \\
	Guanajuato \\
	Mexico}
\email{quiroga@cimat.mx}
\author{A.~S\'anchez-Nungaray}
\address{Facultad de Matem\'aticas \\
	Universidad Veracruzana \\
	Veracruz \\
	Mexico}
\email{armsanchez@uv.mx}

\subjclass{47B35, 32A36, 22E46, 32M15}
\keywords{Toeplitz operators, Projective Space}
\thanks{Research supported by SNI and Conacyt Grants.}

\begin{abstract}
	We consider separately radial (with corresponding group $\T^n$) and radial (with corresponding group $\U(n))$ symbols on the projective space $\PnC$, as well as the associated Toeplitz operators on the weighted Bergman spaces. It is known that the $C^*$-algebras generated by each family of such Toeplitz operators are commutative (see \cite{QS-Projective}). We present a new representation theoretic proof of such commutativity. Our method is easier and more enlightening as it shows that the commutativity of the $C^*$-algebras is a consequence of the existence of multiplicity-free representations. Furthermore, our method shows how to extend the current formulas for the spectra of the corresponding Toeplitz operators to any closed group lying between $\T^n$ and $\U(n)$. 
\end{abstract}

\maketitle

\section{Introduction}\label{sec:Introduction}
The study of families of Toeplitz operators with special properties is already prominent in operator theory. This is particularly so given the abundance of commutative algebras generated by Toeplitz operators on complex spaces. With this respect, the results from \cite{DOQJFA} have special relevance as they prove the existence of nontrivial commutative $C^*$-algebras generated by Toeplitz operators on every weighted Bergman space over every irreducible bounded symmetric domain. The first fundamental ingredient to prove the latter is one which has become very important for Toeplitz operators in the last years: group theory. The origins of this fact include a classification of commutative $C^*$-algebras generated by Toeplitz operators on the unit disk (see \cite{GQV}) and the first examples of large families of commuting Toeplitz operators on the unit ball (see \cite{QVBall1} and \cite{QVBall2}). The second fundamental ingredient for the findings presented in \cite{DOQJFA} is the representation theory of semisimple Lie groups and its study of the holomorphic discrete series for simple Lie groups of Hermitian type.

On the other hand, the existence of commutative $C^*$-algebras of Toeplitz operators on Bergman spaces over complex projective spaces has been established in \cite{QS-Projective}. The Toeplitz operators in such previous work are those whose symbols are invariant under the action of a maximal torus of isometries. These symbols are called separately radial (see Section~\ref{sec:Toeplitz_separately_radial} for further details). The use of group theory in the results obtained in \cite{QS-Projective} is already patent in the definition of the separately radial symbols in terms of group actions. Nevertheless, the proof found in \cite{QS-Projective} for the commutativity of the $C^*$-algebras generated by Toeplitz operators with separately radial symbols is a rather straightforward computation. In other words, no representation theory is involved. 

The first main goal of this work is to explore the relationship between the representation theory associated to the projective space $\PnC$ and the commutativity of the $C^*$-algebra $\mT^{(m)}(\mA^{\T^n})$ generated by Toeplitz operators with separately radial symbols. We observe that $\T^n$ stands for the $n$-dimensional torus acting componentwise on inhomogeneous coordinates of $\PnC$; we refer to Section~\ref{sec:Toeplitz_separately_radial} for further details and for the moment we note that $\mA^K$ denotes the space of symbols invariant under a group $K$. Our first goal is achieved in Section~\ref{sec:Toeplitz_separately_radial}, where we present in Theorem~\ref{thm:Toep_Sep_Radial} a purely representation theoretic proof of the commutativity of $\mT^{(m)}(\mA^{\T^n})$. A comparison with the corresponding proof found in \cite{QS-Projective} clearly shows the efficiency of our new representation theoretic proof. Furthermore, our new proof is more enlightening as it shows the reason for the commutativity of these $C^*$-algebras: the presence of a multiplicity-free action of the group $\T^n$ (see Section~\ref{sec:Invariant_symbols} for further details). Also, our method continues to provide the simultaneous diagonalization that allows to analyze the spectra of the Toeplitz operators involved.

An important case considered in \cite{DOQJFA} is the set of symbols that are invariant under the maximal compact subgroups of isometries of an irreducible bounded symmetric domain. For the projective space $\PnC$, the group of isometries is given by an action of the projective unitary group $\PU(n+1)$ for which the subgroup of isometries fixing a given point is (up to conjugacy) realized by the subgroup $\U(n)$. The second main goal of this paper is to consider symbols invariant under the action of $\U(n)$ (radial symbols) and establish the commutativity of the $C^*$-algebras $\mT^{(m)}(\mA^{\U(n)})$ generated by Toeplitz operators with such radial symbols. Again, our method is representation theoretic and far easier than a corresponding straightforward proof (compare with the methods of \cite{GKV}), as well as more enlightening since a multiplicity-free representation is still involved but now for $\U(n)$. The main result in this case is Theorem~\ref{thm:Toep_Radial}.

On the other hand, if we let $K$ be a closed group such that $\T^n \subset K \subset \U(n)$, then for the above notation we have $\mA^{\U(n)} \subset \mA^K \subset \mA^{\T^n}$. This implies the corresponding inclusions $\mT^{(m)}(\mA^{\U(n)}) \subset \mT^{(m)}(\mA^K) \subset \mT^{(m)}(\mA^{\T^n})$. In particular, the commutativity (however it is proved) of the $C^*$-algebra generated by Toeplitz operators with separately radial symbols implies the commutativity of the corresponding $C^*$-algebra for the group $K$. Nevertheless, our representation theoretic proof has the advantage of providing a general method to describe the spectra of the Toeplitz operators in $\mT^{(m)}(\mA^K)$ as long as the group $K$ is explicitly given: this is the content of Theorem~\ref{thm:Toep_spec_K-inv}.

\section{Bergman spaces and Toeplitz operators on the projective space}\label{sec:BergmanToeplitz}
We denote by $\PnC$ the $n$-dimensional complex projective space which consists of the classes $[w] = \C w \setminus \{0\}$ where $w \in \C^{n+1}$. For every $j = 0, \dots, n$ we consider the inhomogeneous coordinates $\varphi_j : U_j \rightarrow \C^n$ given by
\[
	[w] \mapsto (z_1, \dots, z_n) = \frac{1}{w_j}\left(w_0 \dots, \widehat{w}_j, \dots, w_n\right), 
\]
where $U_j$ is the open set of the classes $[w]$ such that $w_j \not= 0$. 

We will denote by $\omega$ the canonical symplectic form on $\PnC$ which yields the Fubini-Study geometric structure (see \cite{QS-Projective} for further details). And $\Omega$ will denote the corresponding volume form on $\PnC$.

If we consider $\T$ as the subgroup of $\U(n+1)$ that consists of multiples of the identity, then we have a natural quotient map of Lie groups
\[
	\U(n+1) \rightarrow \U(n+1)/\T = \PU(n+1)
\]
where $\PU(n+1)$ is the projective unitary group. The image of a matrix $A \in \U(n+1)$ will be denoted by $[A]$.
Recall that $\PU(n+1)$ is precisely the group of biholomorphic isometries of $\PnC$ through the following action
\begin{align*}
	\PU(n+1) \times \PnC &\rightarrow \PnC \\
		([A], [z]) &\mapsto [Az].
\end{align*}
The isotropy subgroup for this action corresponding to the point $[(1,0,\dots, 0)^t]$ is given by the subgroup
\[
	\left\{
		\begin{bmatrix}
		t & 0 \\
		0 & A
		\end{bmatrix} \sVert[3] t \in \U(1), A \in \U(n)
	\right\}
	=\left\{
		\begin{bmatrix}
			1 & 0 \\
			0 & A
		\end{bmatrix} \sVert[3] A \in \U(n)
	\right\}.
\]
Since this is canonically isomorphic to $\U(n)$, we will use the latter to denote both the isotropy subgroup and the corresponding group of $n\times n$ matrices. On the other hand, by the corresponding quotient map, the group $\U(n+1)$ acts on $\PnC$ biholomorphically with isotropy subgroup for the point $[(1,0, \dots, 0)^t]$ given by
\[
	\U(1) \times \U(n) =
	\left\{
		\begin{pmatrix}
			t & 0 \\
			0 & A
		\end{pmatrix} \sVert[3] t \in \U(1), A \in \U(n)
	\right\}.
\]
In particular, we have
\[
	\PnC \cong \PU(n+1)/\U(n) \cong \U(n+1)/(\U(1) \times \U(n)).
\]

For every $m \in \Z$, we consider the character $\chi_m : \U(1) \times \U(n) \rightarrow \T$ given by $\chi_m(t) = t^m$. This yields the associated line bundle over $\PnC$ defined by
\[
	L_m = \U(n+1) \times_{\chi_m} \C \rightarrow \PnC.
\]
As usual, the elements of this bundle are the classes $[A, \lambda]$ for the equivalence relation on $\U(n+1) \times \C$ given by
\[
	(A_1, \lambda_1) \sim (A_2, \lambda_2) \iff (A_2, \lambda_2) = (A_1 B, \chi_m(B)^{-1}\lambda_1) \text{ for some } B \in \U(1) \times \U(n).
\]
By construction, these line bundles admit $\U(n+1)$-actions that preserve Hermitian metrics inherited from the metric on the standard fiber $\C$.

We are interested in the hyperplane line bundle $H = L_{-1}$ and its tensor powers $H^m = L_{-m}$ for $m \in \N_0$. These are precisely the line bundles over $\PnC$ that admit holomorphic sections. For every $m \in \N_0$, the space of holomorphic sections of $H^m$ is called the weighted Bergman space $\mA^2_m(\PnC)$ with weight $m$.

If we denote by $h^m$ the Hermitian metric of $H^m$, then we can define $L^2(\PnC,H^m)$ as the Hilbert space of measurable sections $f$ of $H^m$ that satisfy
\[
\int_{\PnC} h^m(f,f) \Omega < \infty,
\]
and endowed with the corresponding inner product. It follows that $\mA^2_m(\PnC)$ is a finite dimensional and so closed subspace of the Hilbert space $L^2(\PnC,H^m)$.

The line bundles $H^m$ can be trivialized on every open set $U_j$ as above. For example, on the open set $U_0$ the holomorphic map
\begin{align*}
	H^m|_{U_0} &\rightarrow \C^n \times \C \\
	[A,\lambda] &\mapsto \left([Ae_1],\frac{\lambda}{w_0^m}\right)
\end{align*}
yields a trivialization, where $Ae_1 = (w_0, \dots, w_n)$. It is well known that, with respect to this trivialization, the Bergman space $\mA^2_m(\PnC)$ corresponds precisely to the space $\mP_m(\C^n)$ of polynomials with degree at most $m$ on $\C^n$. 

We will use the trivialization $H^m|_{U_0} \cong \C^n\times \C$ and the corresponding identifications for the spaces of sections to perform our computations using inhomogeneous coordinates. In particular, for every $m \in \N_0$ let us define the measure on $\C^n$ given by
\[
	\dif\nu_m(z) = \frac{(n+m)!}{\pi^n m!} \frac{\dif z}{(1+|z_1|^2 + \dots + |z_n|^2)^{n+m+1}}.
\]
This yields the Hilbert space $L^2(\C^n, \nu_m)$ which clearly contains $\mP_m(\C^n)$ as a finite dimensional and so closed subspace. It follows from \cite{QS-Projective} that for every $m \in \N_0$ there is an isometry $L^2(\PnC, H^m) \rightarrow L^2(\C^n, \nu_m)$ obtained from the trivialization of $H^m|_{U_0}$ described above that maps $\mA^2_m(\PnC)$ onto $\mP_m(\C^n)$. Furthermore, the Bergman projection $B_m : L^2(\C^n, \nu_m) \rightarrow \mP_m(\C^n)$ is given by
\[
	B_m(f)(z) = \int_{\C^n} f(w) K_m(z,w) \dif\nu_m(w),
\]			
where $K_m(z,w) = (1 + z_1 \overline{w}_1 + \dots + z_n \overline{w}_n)^m$.

In particular, in the realization of the Bergman spaces and Bergman projections described above, the Toeplitz operator with symbol $a \in L^\infty(\C^n)$ is given by
\begin{align*}
	T^{(m)}_a : \mP_m(\C^n) &\rightarrow \mP_m(\C^n) \\
		T^{(m)}_a(f)(z) &= \frac{(n+m)!}{\pi^n m!} \int_{\C^n} \frac{a(w)f(w) (1 + z_1 \overline{w}_1 + \dots + z_n \overline{w}_n)^m \dif w}{(1+|w_1|^2 + \dots + |w_n|^2)^{n+m+1}}.
\end{align*}

Next we observe that the natural isometric action of $\U(n+1)$ on the line bundle $H^m$  induces a unitary representation on the Hilbert space $L^2(\PnC, H^m)$. This representation leaves invariant the corresponding Bergman space $\mA^2_m(\PnC)$ which, by Bott-Borel-Weil theorem, yields an irreducible unitary representation of $\U(n+1)$. 

We note that, for the projection map $\U(n+1) \rightarrow \PU(n+1)$, the subgroup $\U(n)$ of $\PU(n+1)$ is already the image of the subgroup $\{1\} \times \U(n) \subset \U(1) \times \U(n) \subset \U(n+1)$. This allows to obtain a representation of the subgroup $\U(n)$ of $\PU(n+1)$ on the space of sections of $H^m$. Furthermore, for the trivialization of $H^m|_{U_0}$ we conclude that the corresponding action is given by
\begin{equation}\label{eq:U(n)action}
(A \cdot f)(z) = f(A^{-1}z),
\end{equation}
for every $A \in \U(n)$, $f : \C^n \rightarrow \C$ and $z \in \C^n$. As a particular case, the expression \eqref{eq:U(n)action} yields a unitary representation of $\U(n)$ on $L^2(\C^n, \nu_m)$ (corresponding to the representation on  $L^2(\PnC, H^m)$) that preserves $\mP_m(\C^n)$ (corresponding to the Bergman space $\mA^2_m(\PnC)$). 

In the rest of this work we will make use of such realizations of these unitary representations of $\U(n)$. And we will denote by $\pi_m$ the unitary representation of $\U(n)$ on the Bergman space $\mA^2_m(\PnC) \cong \mP_m(\C^n)$.

\section{Invariant symbols, Toeplitz operators and intertwining maps}\label{sec:Invariant_symbols}
Let us fix a closed subgroup $K$ of $\U(n)$ and consider the vector subspace $\mA^K$ of $L^\infty(\PnC)$ of $K$-invariant symbols. More precisely, we have $a \in \mA^K$ if and only if 
\[
	a(kz) = a(z) 
\]
for all $k \in K$ and for a.e.~$z \in \C^n$. Correspondingly, for every $m \in \N_0$ we will denote by $\mT^{(m)}(\mA^K)$ the $C^*$-algebra generated by the Toeplitz operators on $\mA^2_m(\PnC)$ with symbols in $\mA^K$.

Since $K$ is a subgroup of $\U(n)$, the unitary representation $\pi_m$ of $\U(n)$ on $\mA^2_m(\PnC)$ restricted to $K$ is itself a unitary representation which we will denote by $\pi_m|_K$. The following result is an easy consequence of the definitions (see Corollary~3.3 from \cite{DOQJFA}).

\begin{proposition}\label{prop:InvSymbols_IntertwiningToeplitz}
	Let $K$ be a closed subgroup of $\U(n)$ and $a \in L^\infty(\PnC)$ be a bounded symbol. If the symbol $a$ belongs to $\mA^K$, then for every $m \in \N_0$ the Toeplitz operator $T^{(m)}_a$ intertwines the restriction $\pi_m|_K$, in other words, we have
	\[
		\pi_m(k) \circ T^{(m)}_a = T^{(m)}_a \circ \pi_m(k)
	\]
	for all $k \in K$.
\end{proposition}

For $K$ a closed subgroup of $\U(n)$ as above, we denote by $\End_K(\mA^2_m(\PnC))$ the space of $\pi_m|_K$-intertwining linear maps of $\mA^2_m(\PnC)$, i.e.~the linear maps $T : \mA^2_m(\PnC) \rightarrow \mA^2_m(\PnC)$ that satisfy
\[
	\pi_m(k) \circ T = T \circ \pi_m(k),
\]
for all $k \in K$.

With the previous notation, Proposition~\ref{prop:InvSymbols_IntertwiningToeplitz} establishes that 
\[
	\mT^{(m)}(\mA^K) \subset \End_K(\mA^2_m(\PnC)),
\]
for every closed subgroup $K \subset \U(n)$. On the other hand, by definition we clearly have 
\[
	\mT^{(m)}(\mA) \subset \End(\mA^2_m(\PnC)).
\]
We will see that these inclusions are in fact identities. 

The following result can be obtained by an easy adaption of the arguments found in \cite{Englis}. See also section~5 from \cite{MOS} for an alternative proof.

\begin{theorem}\label{thm:Toeplitz=Linear}
	Every linear map $\mA^2_m(\PnC) \rightarrow \mA^2_m(\PnC)$ can be realized as a Toeplitz operator with bounded symbol. In other words, we have
	\[
		\End(\mA^2_m(\PnC)) = \mT^{(m)}(L^\infty(\PnC)).
	\]
\end{theorem}

For the next result we will follow the arguments found in \cite{DOQJFA}. For this we introduce some notation. 

Let $K$ be a closed subgroup of $\U(n)$. For every $a \in L^\infty(\PnC)$ let us denote
\[
	\widetilde{a}(z) = \int_K a(k^{-1}z) \dif k,
\]
where the integral is computed with respect to the Haar measure of $K$. Note that $\widetilde{a} \in L^\infty(\PnC)$ by the compactness of $K$. Furthermore, the invariance of the Haar measure clearly implies that $\widetilde{a} \in \mA^K$ for every $a \in L^\infty(\PnC)$. We also have $\widetilde{\widetilde{a}} = \widetilde{a}$, for every $a \in L^\infty(\PnC)$. In particular, the map $a \mapsto \widetilde{a}$ defines a linear projection $L^\infty(\PnC) \rightarrow \mA^K$.

For $K$ as above, a similar construction can be performed for linear maps on $\mA^2_m(\PnC)$. For a linear map $T : \mA^2_m(\PnC) \rightarrow \mA^2_m(\PnC)$ we define
\[
	\widetilde{T} = \int_K \pi_m(k) \circ T \circ \pi_m(k)^{-1} \dif k.
\]
Hence, it is easily seen that $\widetilde{T} \in \End_K(\mA^2_m(\PnC))$ and that the map $T \mapsto \widetilde{T}$ defines a linear projection $\End(\mA^2_m(\PnC)) \rightarrow \End_K(\mA^2_m(\PnC))$.

A straightforward computation proves (see \cite{DOQJFA}) that for every $a \in L^\infty(\PnC)$ and $m \in \N_0$ we have 
	\[
		\widetilde{T^{(m)}_a} = T^{(m)}_{\widetilde{a}}.
	\]

As a consequence we obtain the following result.

\begin{theorem}\label{thm:KToeplitz=KLinear}
	Let $K$ be closed subgroup of $\U(n)$. Every $K$-intertwining linear map $\mA^2_m(\PnC) \rightarrow \mA^2_m(\PnC)$ can be realized as a Toeplitz operator with a $K$-invariant bounded symbol. In other words, we have
	\[
		\End_K(\mA^2_m(\PnC)) = \mT^{(m)}(\mA^K).
	\]
\end{theorem}
\proof
	Let $T \in \End_K(\mA^2(\PnC))$ be given. By Theorem~\ref{thm:Toeplitz=Linear} there exists $a \in L^\infty(\PnC)$ such that $T = T^{(m)}_a$. And so we have
	\[
		T = \widetilde{T} = \widetilde{T^{(m)}_a} = T^{(m)}_{\widetilde{a}}.
	\]
	Since $\widetilde{a} \in \mA^K$ the result follows.
\endproof

The previous result reduces the study of the $C^*$-algebra $\mT^{(m)}(\mA^K)$ to the study of the representation of $K$ on $\mA^2_m(\PnC)$.

Let us fix a closed subgroup $K$ of $\U(n)$. Since $K$ is compact and $\mA^2_m(\PnC)$ is finite dimensional there is a collection $U_1, \dots, U_l$ of mutually inequivalent irreducible $K$-modules and a $K$-invariant direct sum decomposition
\begin{equation}\label{eq:decomp_irred}
	\mA^2_m(\PnC) = \bigoplus_{j=1}^{l} V_j
\end{equation}
so that the subspace $V_j$ is a direct sum of irreducible $K$-submodules isomorphic to $U_j$. This is called the isotypic decomposition of the $K$-module $\mA^2_m(\PnC)$. We recall that the number of linearly independent summands in $V_j$ isomorphic to $U_j$ is called the multiplicity (of the latter in the former).

With respect to the isotypic decomposition \eqref{eq:decomp_irred}, the fact that the $K$-modules $U_j$ are mutually inequivalent implies that
\begin{equation}\label{eq:EndK_isotypic_comp}
	\End_K(\mA^2_m(\PnC)) = \bigoplus_{j=1}^{l} \End_K(V_j),
\end{equation}
as complex algebras. Hence, to understand the structure of $\End_K(\mA^2_m(\PnC))$ it is enough to consider a $K$-module $V$ given by a direct sum of copies of an irreducible $K$-module and describe $\End_K(V)$. And we recall that $\End_K(V) \cong M_r(\C)$, where $r$ is the multiplicity of $U$ in $V$.

From the previous discussion we conclude the following result. We recall that the representation of $K$ on $\mA^2_m(\PnC)$ is called multiplicity-free if the isotypic decomposition~\eqref{eq:decomp_irred} has multiplicity $1$ on each summand.

\begin{theorem}\label{thm:Toeplitz_K-inv}
	Let $K$ be a closed subgroup of $\U(n)$. Then, we have
	\[
		\mT^{(m)}(\mA^K) \cong \bigoplus_{j=1}^{l} M_{r_j}(\C),
	\]
	as $C^*$-algebras, where $r_j$ is the multiplicity of $U_j$ in $V_j$ in the isotypic decomposition~\eqref{eq:decomp_irred}.
	
	In particular, for every $m \in \N_0$ the following conditions are equivalent.
	\begin{enumerate}
		\item The $C^*$-algebra $\mT^{(m)}(\mA^K)$ is commutative.
		\item The representation of $K$ on $\mA^2_m(\PnC)$ is multiplicity-free.
	\end{enumerate}
\end{theorem}

\section{Toeplitz operators with separately radial symbols}\label{sec:Toeplitz_separately_radial}
Here we apply the previous discussion to the case $K = \T^n$ the subgroup of $\U(n)$ consisting of the diagonal matrices. The elements of $\mA^{\T^n}$ are called separately radial symbols. The results in this section may be compared with those from \cite{QS-Projective}, particularly Section~5 and Theorem~5.5 therein.

To apply Theorem~\ref{thm:Toeplitz_K-inv} we obtain the isotypic decomposition of the representation of $\T^n$ on $\mA^2_m(\PnC)$. We recall from Section~\ref{sec:BergmanToeplitz} the identification of $\mA^2_m(\PnC)$ with the space $\mP_m(\C^n)$ of complex polynomials on $\C^n$ with dimension at most $m$.

We will freely use the multi-index notation and we will also consider the set
\[
	J_n(m) = \{p \in \N_0^n \sVert |p| \leq m\},
\]
where $m \in \N_0^m$.

Hence, the set of monomials $z^p$, where $p \in J_n(m)$, is a basis for $\mP_m(\C^n)$. Furthermore, we have
\[
	\pi_m|_{\T^n}(t) z^p = t^{-p} z^p,
\]
for every $t \in \T^n$ and $p \in J_n(m)$. It follows that the decomposition into irreducible submodules of the representation of $\T^n$ on $\mP_m(\C^n)$ is given by
\begin{equation}\label{eq:iso_decomp_Tn}
	\mP_m(\C^n) = \bigoplus_{p \in \N_0, |p| \leq m} \C z^p.
\end{equation}
Moreover, for $p \in \N_0$ with $|p|\leq m$, the submodule $\C z^p$ has character $\chi_p$ given by
\begin{align*}
	\chi_p : \T^n &\rightarrow \T \\
		t &\mapsto t^{-p}.
\end{align*}
Since all such characters are distinct, it follows that \eqref{eq:iso_decomp_Tn} is in fact the isotypic decomposition, which is then multiplicity-free. 

On the other hand, an easy computation shows that (see \cite{QS-Projective})
\[
	\left\{e_p(z) = \left(\frac{m!}{p! (m-|p|)!}\right)^{\frac{1}{2}} z^p \sVert[4] p \in J_n(m) \right\}
\]
is an orthonormal basis. In particular, the map given by
\begin{align*}
	R : \mA^2_m(\PnC) &\rightarrow \ell^2(J_n(m)) \\
	f &\mapsto \left(\left<f,e_p\right>_m\right)_{p \in J_n(m)},
\end{align*}
is an isometry.

We now present an elementary proof of Theorem~5.5 from \cite{QS-Projective}.

\begin{theorem}\label{thm:Toep_Sep_Radial}
	For every $m \in \N_0$ the $C^*$-algebra $\mT^{(m)}(\mA^{\T^n})$ of Toeplitz operators with $\T^n$-invariant symbols is commutative. Furthermore, for every $a \in \mA^{\T^n}$ we have
	\[
		R T^{(m)}_a R^* = \gamma_{a,m} I
	\]
	where the function $\gamma_{a,m} : J_n(m) \rightarrow \C$ is given by
	\begin{align*}
		\gamma_{a,m}(p) &= \left<T^{(m)}_a(e_p), e_p\right>_m = \left<a e_p, e_p\right>_m \\
		&= \frac{2^n m!}{p!(m-|p|)!}\int_{R^n_+} 
		\frac{a(r_1, \dots, r_n) r^{2p_1+1} \cdots r^{2p_n+1}\dif r_1 \cdots \dif r_n}{(1+r_1^2+\dots+r_n^2)^{n+m+1}} \\
		&= \frac{m!}{p!(m-|p|)!}\int_{R^n_+} 
		\frac{a(\sqrt{r_1}, \dots, \sqrt{r_n}) r^{p_1+1} \cdots r^{p_n+1}\dif r_1 \cdots \dif r_n}{(1+r_1+\dots+r_n)^{n+m+1}}.
	\end{align*}
\end{theorem}
\proof
	That $R$ is an isometry follows from the fact that the monomials $e_p$ form an orthonormal basis.
	
	On the other hand, using the isotypic decomposition \eqref{eq:iso_decomp_Tn}, Theorem~\ref{thm:Toeplitz_K-inv} implies that $\mT^{(m)}(\mA^{\T^n})$ is commutative. Furthermore, it follows that for every $a \in \mA^{\T^n}$ and $p \in J_n(m)$ we have
	\[
		T^{(m)}_a(e_p) = \gamma_{a,m}(p) e_p
	\]
	for some complex number $\gamma_{a,m}(p)$, thus also showing that
	\[
		\gamma_{a,m}(p) = \left<T^{(m)}_a(e_p), e_p\right>_m = \left<a e_p, e_p\right>_m.
	\]  
	If we denote by $\widetilde{e}_p$ ($p \in J_n(m)$) the elements of the canonical basis of $\ell^2(J_n(m))$, then $R(e_p) = \widetilde{e}_p$ for every $p$ and so we have
	\[
		R T^{(m)}_a R^*(\widetilde{e}_p) = R T^{(m)}_a(e_p) 
			= R(\gamma_{a,m}(p)e_p) 
			= \gamma_{a,m}(p)\widetilde{e}_p,
	\]
	which shows that $R T^{(m)}_a R^* = \gamma_{a,m} I$.
	
	Finally, to write down the function $\gamma_{a,m}$ explicitly we use
	\[
		\gamma_{a,m}(p) = \left< a e_p, e_p \right>_m 
			= \frac{(n+m)!}{\pi^n p!(m - |p|)!}\int_{\C^n} \frac{a(z)|z^p|^2 \dif z}{(1 + |z|^2)^{n+m+1}}
	\]
	and apply polar coordinates to each complex coordinate to obtain the first integral expression. This uses the fact that $a(z_1, \dots, z_n) = a(|z_1|, \dots, |z_n|)$. The second integral expression is obtained with the change of coordinates $r \mapsto r^2$.
\endproof

\section{Toeplitz operators with radial symbols}\label{sec:Toeplitz_radial}
We now consider the case $K = \U(n)$, the group of isometries of $\PnC$ that fix the point $[(1,0,\dots, 0)^t]$. The elements of $\mA^{\U(n)}$ are called radial symbols. The results in this section are dual to those considered in \cite{GKV} and \cite{Q-Grudsky60}.

To understand this case, we need to obtain the isotypic decomposition \eqref{eq:decomp_irred} for $K = \U(n)$. By the identification of $\mA^2_m(\PnC)$ with $\mP_m(\C^n)$ the following result provides such decomposition.

\begin{proposition}\label{prop:Iso_decomp_radial}
	For every $k \in \N_0$, denote by $\mP^k(\C^n)$ the space of homogeneous polynomials with degree $k$ on $\C^n$.
	Then, the isotypic decomposition of the representation of $\U(n)$ on $\mP_m(\C^n)$ is given by
	\[
		\mP_m(\C^n) = \bigoplus_{k=0}^m \mP^k(\C^n),
	\]
	where the summands $\mP^k(\C^n)$ are irreducible and inequivalent as $\U(n)$-modules. Hence, the representation of $\U(n)$ on $\mA^2_m(\C^n)$ is multiplicity-free for every $m \in \N_0$.
\end{proposition}
\proof
	The direct sum above holds trivially and it is a decomposition into $\U(n)$-submodules since the representation of $\U(n)$ is linear on the variable $z \in \C^n$. Hence, it is enough to show that the spaces $\mP^k(\C^n)$ are irreducible and inequivalent over $\U(n)$.
	
	First we note that the $\U(n)$ representation on $\mP^k(\C^n)$ is the restriction of the rational representation 
	\begin{align*}
		\GL(n,\C) \times \mP^k(\C^n) &\rightarrow \mP^k(\C^n) \\
			(A, f(z)) &\mapsto f(A^{-1}z).
	\end{align*}
	Furthermore, by classical invariant theory (see \cite{GoodmanWallach}) it is known that $\mP^k(\C^n)$ is an irreducible $\GL(n,\C)$-module for every $k \geq 0$. On the other hand, $\U(n)$ is a real form and so a Zariski dense subgroup of $\GL(n,\C)$. In particular, any subspace of $\mP^k(\C^n)$ is invariant under $\GL(n,\C)$ if only if it is invariant under $\U(n)$. As a consequence, $\mP^k(\C^n)$ is an irreducible module over $\U(n)$ as well. Finally, since the dimensions of the spaces $\mP^k(\C^n)$ are all different these spaces are necessarily inequivalent over $\U(n)$ by their irreducibility.
\endproof

As a consequence of the previous results, if $T \in \End_{\U(n)}(\mA^2_m(\PnC))$, then we have
\[
	T(\mP^k(\C^n)) \subset \mP^k(\C^n)
\] 
for every $k = 0, \dots, m$, with $T$ acting as a multiple of the identity $c_k I$ on every such space. If for each $k$ we let $u_k \in \mP^k(\C^n)$ be a unitary vector, then we can compute the scalar $c_k$ by
\[
	c_k = \left<T u_k, u_k \right>_m. 
\]

For our purposes it will be useful to consider for every $k = 0, \dots, m$ the homogeneous polynomial
\[
	f_k(z) = \sum_{\substack{p \in \N_0^n\\|p| = k}} \sqrt{\binom{k}{p}} z^p
\]
where the multinomial coefficients are defined by
\[
	\binom{k}{p} = \frac{k!}{p_1!\cdots p_n!},
\]
for every $p \in \N_0^n$ such that $|p| = k$.

\begin{lemma}\label{lem:norm_poly}
	For every $k,m \in \N_0$ with $k \leq m$, we have
	\[
		\|f_k\|_m^2 = \frac{2(n+m)!}{m!(n-1)!}\int_{0}^{+\infty}\frac{r^{2n+2k-1}}{(1+r^2)^{n+m+1}}\dif r.
	\]
\end{lemma}
\proof
	For every $k \leq m$ we compute
	\begin{align*}
		\left<f_k, f_k\right>_m &= \sum_{\substack{p,q \in \N_0^n\\|p|=|q|=k}} \sqrt{\binom{k}{p}} \sqrt{\binom{k}{q}} \frac{(n+m)!}{\pi^n m!} \int_{\C^n} \frac{z^p \overline{z}^q \dif z}{(1+|z|^2)^{n+m+1}} \\
	\intertext{by the orthonormality of the monomials used in Section~\ref{sec:Toeplitz_separately_radial} we have}
		&= \frac{(n+m)!}{\pi^n m!} \int_{\C^n} \sum_{\substack{p \in \N_0^n\\|p|=k}} \binom{k}{p} |z_1|^{2p_1} \cdots |z_n|^{2p_n} \frac{\dif z}{(1+|z|^2)^{n+m+1}} \\
		&= \frac{(n+m)!}{\pi^n m!} \int_{\C^n} (|z_1|^2 + \dots + |z_n|^2)^k \frac{\dif z}{(1+|z|^2)^{n+m+1}} \\
	\intertext{introducing spherical coordinates with $\dif \sigma$ the normalized volume of the unit sphere $\mathbb{S}^n$ in $\C^n$ we obtain}
		&= \frac{(n+m)!}{\pi^n m!} \int_{0}^{+\infty} \int_{\mathbb{S}^n} \frac{r^{2k}}{(1+r^2)^{n+m+1}} 
		\frac{2\pi^n}{(n-1)!} r^{2n-1}\dif r \dif \sigma \\
		&= \frac{2 (n+m)!}{m!(n-1)!} \int_{0}^{+\infty} \frac{r^{2n + 2k - 1}}{(1+r^2)^{n+m+1}}  \dif r.
	\end{align*}
\endproof

Note that the set
\[
	\left( u_k = \frac{1}{\|f_k\|_m} f_k \right)_{k = 0, \dots, m}
\]
is orthonormal in $\mA^2_m(\PnC)$ with exactly one vector chosen from each one of the summands in the decomposition from Proposition~\ref{prop:Iso_decomp_radial}. And Lemma~\ref{lem:norm_poly} provides an explicit formula for the norms $\|f_k\|_m$. 

\begin{lemma}\label{lem:ortho_with_a}
	If $m \in \N_0$ and $a \in \mA^{\U(n)}$ are given, then
	\[
		\left< a e_p, e_q\right>_m = 0
	\]
	for any two different elements $p,q \in J_n(m)$.
\end{lemma}
\proof
	Let $p \not= q$ as above be given. Hence, with the character notation used in Section~\ref{sec:Toeplitz_separately_radial} there exists $t \in \T^n$ such that
	\[
		\chi_p(t) \not= \chi_q(t).
	\]
	Consider the unitary transformation given by $z \mapsto (t_1 z_1, \dots, t_n z_n)$. Since the measure $a(z) \dif \nu_m(z)$ is invariant under unitary transformations we conclude that
	\begin{align*}
		\left< a e_p, e_q\right>_m & = \int_{\C^n} e_p(z) \overline{e_q(z)} a(z) \dif \nu_m(z) \\
			&= \chi_p(t) \overline{\chi_q(t)} \int_{\C^n} e_p(z) \overline{e_q(z)} a(z) \dif \nu_m(z) \\
			&= \chi_p(t) \overline{\chi_q(t)} \left< a e_p, e_q\right>_m,
	\end{align*}
	from which our result follows.
\endproof

\begin{theorem}\label{thm:Toep_Radial}
	For every $m \in \N_0$ the $C^*$-algebra $\mT^{(m)}(\mA^{\U(n)})$ of Toeplitz operators with $\U(n)$-invariant symbols is commutative. Furthermore, for every $a \in \mA^{\U(n)}$ we have
	\[
		R T^{(m)}_a R^* = \gamma_{a,m} I,
	\]
	where the complex function $\gamma_{a,m} : J_n(m)  \rightarrow \C$ satisfies
	\[
		\gamma_{a,m}(p) = \gamma_{a,m}(q),
	\]
	whenever $|p|=|q|$. The induced complex function $\widehat{\gamma}_{a,m}$ defined by
	\[
		\widehat{\gamma}_{a,m}(|p|) = \gamma_{a,m}(p)
	\]
	for $p \in J_n(m)$ can be computed by
	\begin{align*}
		\widehat{\gamma}_{a,m}(k) &= \left<T^{(m)}_a u_k, u_k \right>_m = \left< a u_k, u_k \right>_m \\
			&= \frac{\displaystyle\int_{0}^{+\infty}\frac{a(r)r^{2n+2k-1}\dif r}{(1+r^2)^{n+m+1}}}{\displaystyle\int_{0}^{+\infty}\frac{r^{2n+2k-1}\dif r}{(1+r^2)^{n+m+1}}} 
			= \frac{(n+m)!}{ (n+k-1)! (m-k)!}     \displaystyle\int_{0}^{+\infty}\frac{a(\sqrt{r})r^{n+k-1}\dif r}{(1+r)^{n+m+1}}  ,
	\end{align*}
	for every $k \leq m$.
\end{theorem}
\proof
	The commutativity of $\mT^{(m)}(\mA^{\U(n)})$ follows from Theorem~\ref{thm:Toeplitz_K-inv} and Proposition~\ref{prop:Iso_decomp_radial}. Let us now choose $a \in \mA^{\U(n)}$.
	
	Since the irreducible summands in the isotypic decomposition have orthogonal bases given by subsets of $\left(e_p\right)_{p \in J_n(m)}$ the existence of $\gamma_{a,m}$ follows with the same proof as in Theorem~\ref{thm:Toep_Sep_Radial}. 
	
	On the other hand, as noted above in this section we have
	\[
		T^{(m)}_a u = c u,
	\]
	for every $u \in \mP^k(\C^n)$ where the constant depends only on $k$. Hence, the function $\gamma_{a,m}$ is constant on the collection of $p \in J_n(m)$ such that $|p| = k$ is fixed. In particular, the function $\widehat{\gamma}_{a,m}$ is well defined.
	
	To compute $\widehat{\gamma}_{a,m}$ we note that from the previous remarks we have
	\[
		\widehat{\gamma}_{a,m}(k) = \left<T^{(m)}_a u_k, u_k \right>_m = \left< a u_k, u_k \right>_m
			=\frac{\left< a f_k, f_k \right>_m}{\left< f_k, f_k \right>_m}.
	\]
	Then, the inner product in the denominator is obtained from Lemma~\ref{lem:norm_poly} and the inner product in the numerator can be computed similarly to obtain the first integral expression. The second integral expression is obtained with the change of coordinates $r \mapsto r^2$.
\endproof

\section{Final remarks}\label{sec:final}
In this work we have considered the isometry group $\PU(n+1)$ of the projective space $\PnC$, paying particular attention to the maximal subgroup $\U(n)$. Nevertheless it is natural to ask about Toeplitz operators with symbols invariant under some other closed subgroups not necessarily contained in $\U(n)$. In the picture provided by the group $\U(n+1)$ this corresponds to the maximal compact subgroup $\U(1) \times \U(n)$. And it is also well known that $\U(n+1)$ contains maximal subgroups inequivalent (under conjugation) to $\U(1) \times \U(n)$. Some examples are given by the subgroups $\U(j)\times\U(k)$ where $j+k = n+1$. But we recall that the corresponding quotient $\U(n+1) / (\U(j)\times \U(k))$ is the complex Grassmannian of $j$-planes in $\C^{n+1}$. It seems more reasonable to study this setup in the context of such spaces.

On the other hand, the methods of Sections~\ref{sec:Toeplitz_separately_radial} and \ref{sec:Toeplitz_radial} can be applied to other subgroups of $\U(n)$. More precisely, for a given closed subgroup $K$ of $\U(n)$ and to study the $C^*$-algebra $\mT^{(m)}(\mA^K)$ we first apply Theorem~\ref{thm:Toeplitz_K-inv} to reduce its structure to the representation of $K$ on $\mA^2_m(\PnC) \cong \mP_m(\C^n)$. The important piece of information is the isotypic decomposition described in \eqref{eq:decomp_irred}. When the group $K$ is explicitly given, it is possible to use the character theory of compact groups to write down such decomposition.

The multiplicity-free case is particularly interesting and admits a more direct approach. In such case, let us denote by $q_j = \dim V_j$ the dimension of the irreducible components. Also choose an orthonormal basis $v_1^j, \dots, v_{q_j}^j$ for the space $V_j$. We recall that all such summands are mutually perpendicular since the representation is unitary. Hence, the set
\[
	\{ v_k^j \mid j = 1,\dots,l, \; k=1, \dots, q_j \}
\]
is an orthonormal basis for $\mA^2_m(\PnC)$. As before, we can consider the unitary map
\begin{align*}
	R : \mA^2_m(\PnC) &\rightarrow \ell^2(q_1) \times \dots \times \ell^2(q_l) = \ell^2(|q|) \\
		f &\mapsto \left(\left<f,v_k^j\right>_m \right)_{\substack{j = 1, \dots, l, \\ k = 1, \dots, q_j}}. 
\end{align*}
where $q = (q_1, \dots, q_l)$. Then, the following result provides a generalization of Theorems~\ref{thm:Toep_Sep_Radial} and \ref{thm:Toep_Radial} that can be used to describe the structure of $\mT^{(m)}(\mA^K)$ and the spectra of its elements. For simplicity, we will denote 
\[
	[|q|] = [1:q_1]\dot{\cup} \dots \dot{\cup} [1:q_l]
\]
where $[1:q_j] = \{1, \dots, q_j\}$ and $\dot{\cup}$ denotes disjoint union.

\begin{theorem}\label{thm:Toep_spec_K-inv}
	Let $m \in \N_0$ and $K$ a closed subgroup of $\U(n)$ be given such that the representation of $K$ on $\mA^2_m(\PnC)$ is multiplicity-free so that $\mT^{(m)}(\mA^K)$ is commutative. Then, for every $a \in \mA^K$ and for the notation as above we have
	\[
		R T^{(m)}_a R^* = \gamma_{a,m}I
	\]
	where the function $\gamma_{a,m} : [|q|] \rightarrow \C$ is constant on the interval $[1:q_j]$ for every $j = 1, \dots, l$. Thus inducing a function $\widehat{\gamma}_{a,m} : [1:l] \rightarrow \C$ such that
	\[
		\widehat{\gamma}_{a,m}(j) = \gamma_{a,m}(k)
	\]
	for every $j = 1, \dots, l$ and $k \in [1:q_j]$. Furthermore, if we choose for every $j$ a non-zero element $f_j \in V_j$, then we have
	\[
		\widehat{\gamma}_{a,m}(j) = \frac{\left< T^{(m)}_a f_j, f_j\right>_m}{\left<f_j,f_j\right>_m} 
			= \frac{\left< a f_j, f_j\right>_m}{\left<f_j,f_j\right>_m},
	\]
	for every $j = 1, \dots, l$.
\end{theorem}

The proof of the previous result follows the same arguments used in Sections~\ref{sec:Toeplitz_separately_radial} and \ref{sec:Toeplitz_radial}.


{\small
\begin{thebibliography}{999}

\bibitem{DOQJFA} 
	{\it M.~Dawson, G.~\'Olafsson, R.~Quiroga-Barranco}: 
	Commuting Toeplitz operators on bounded symmetric domains and multiplicity-free restrictions of holomorphic discrete series. 
	J. Funct. Anal. {\bf 268} (2015), no. 7, 1711--1732. Zbl. 1320.47029, MR3315576 

\bibitem{Englis} 
	{\it M.~Engli\v{s}}: 
	Density of algebras generated by Toeplitz operators on Bergman spaces. 
	Arkiv f\"{o}r Matematik {\bf 30} (1992), No. 2, 227--243. Zbl. 0784.46036, MR1289753 

\bibitem{GoodmanWallach} 
	{\it R.~Goodman, N.~R.~Wallach}: 
	Symmetry, representations, and invariants. 
	Graduate Texts in Mathematics, 255. Springer, Dordrecht, 2009. Zbl. 1173.22001, MR2522486 

\bibitem{GKV} 
	{\it S.~Grudsky, A.~Karapetyants, N.~Vasilevski}: 
	Toeplitz operators on the unit ball in $\C^n$ with radial symbols. 
	J. Operator Theory {\bf 49} (2003), no. 2, 325--346. Zbl. 1027.32010, MR1991742 

\bibitem{GQV}  
	{\it S.~Grudsky, R.~Quiroga-Barranco, N.~Vasilevski}: 
	Commutative $C^∗$-algebras of Toeplitz operators and quantization on the unit disk. 
	J. Funct. Anal. {\bf 234} (2006), no. 1, 1--44. Zbl. 1100.47023, MR2214138 

\bibitem{MOS} 
	{\it M.~A.~Morales-Ramos, A.~Sanchez-Nungaray, Josue Ramirez-Ortega}: 
	Toeplitz operators with quasi-separately radial symbols on the complex projective space. 
	Boletin de la Sociedad Matematica Mexicana {\bf 22} (2015), no. 1, 213--227. Zbl. 06562396, MR3473758 

\bibitem{Q-Grudsky60} 
	{\it R.~Quiroga-Barranco}: 
	Separately radial and radial Toeplitz operators on the unit ball and representation theory. 
	To appear in Boletin de la Sociedad Matematica Mexicana.

\bibitem{QS-Projective} 
	{\it R.~Quiroga-Barranco, A.~Sanchez-Nungaray}: 
	Commutative $C^*$-algebras of Toeplitz operators on complex projective space. 
	Integral Equations and Operator Theory {\bf 71} (2011), no. 2, 225--243. Zbl. 1251.47065, MR2838143 

\bibitem{QVBall1}  
	{\it R.~Quiroga-Barranco, N.~Vasilevski}:
	Commutative $C^∗$-algebras of Toeplitz operators on the unit ball. I. Bargmann-type transforms and spectral representations of Toeplitz operators. 
	Integral Equations Operator Theory {\bf 59} (2007), no. 3, 379--419. Zbl. 1144.47024, MR2363015 

\bibitem{QVBall2}  
	{\it R.~Quiroga-Barranco, N.~Vasilevski}:
	Commutative $C^∗$-algebras of Toeplitz operators on the unit ball. II. Geometry of the level sets of symbols. 
	Integral Equations Operator Theory {\bf 60} (2008), no. 1, 89–-132. Zbl. 1144.47025, MR2380317 


\end{thebibliography}
}
\end{document}